\documentclass[12pt]{article}
\usepackage[a4paper, total={7in, 8in}]{geometry}
\usepackage{hyperref}
\usepackage{graphicx}
\usepackage{cite}
\usepackage{comment}
\usepackage{tikz}
\usetikzlibrary{automata,arrows,positioning,calc}
\usepackage{amsthm} 
\usepackage{amsmath}
\usepackage{amssymb}
\usepackage{mathtools}

\newtheorem{theorem}{Theorem}

\newtheorem{prop}[theorem]{Proposition}

\newcommand{\N}{\mathbb{N}}

\newcommand{\veps}{\varepsilon}

\newcommand{\beq}{\begin{equation*}}
\newcommand{\eeq}{\end{equation*}}
\newcommand{\head}[1]
{\markright{\hbox to 0pt{\vtop to 0pt{\hbox{}\vskip 3mm \hrule width
\textwidth \vss} \hss}{\sc #1}}}

\begin{document}

\parindent=0pt
\baselineskip18pt
\parskip6pt
\head{Search problem}

\title{Search for a moving target in a competitive environment\thanks{We would like to thank Aditya Aradhye, Luca Margaritella and Niels Mourmans for their precious comments and suggestions.}}
\date{\today}
\author{Benoit Duvocelle \footnote{Corresponding Author.} \thanks{Department of Quantitative Economics, Maastricht University. Email: b.duvocelle@maastrichtuniversity.nl} \hskip6pt
J\'anos Flesch\thanks{Department of Quantitative Economics, Maastricht University. Email: j.flesch@maastrichtuniversity.nl} \hskip6pt
Hui Min Shi\thanks{Department of Quantitative Economics, Maastricht University. Email: hm.shi@student.maastrichtuniversity.nl} \hskip6pt
Dries Vermeulen\thanks{Department of Quantitative Economics, Maastricht University. Email: d.vermeulen@maastrichtuniversity.nl} \hskip6pt}
\maketitle

\begin{abstract}
We consider a discrete-time dynamic search game in which a number of players compete to find an invisible object that is moving according to a time-varying Markov chain. 
We examine the subgame perfect equilibria of these games. The main result of the paper is that the set of subgame perfect equilibria is exactly the set of greedy strategy profiles, i.e. those strategy profiles in which the players always choose an action that maximizes their probability of immediately finding the object.
We discuss various variations and extensions of the model.\medskip

\noindent
{\bf Keywords:} Game Theory; Search game; Optimal search; Greedy strategy; Subgame perfect equilibrium.\medskip

\end{abstract}

\section{Introduction}

In this paper we consider a dynamic search game, in which an object moves according to a time-varying Markov chain across a finite set of states. The set of competing players can be either finite or infinite. 
At each period, an active player is chosen according to a fixed distribution. The active player can search for the object in one of the possible positions. If the object is found by the active player, this player wins and the game ends. Otherwise, the object moves according to the transition matrix, and the game enters the next period. 
Each player observes the action chosen by his opponents, and the transition probabilities, initial probabilities and probabilities of being the active player are known to all players.
The goal of each player is to find the object and win the game. So each player prefers the winning outcome and is indifferent between the outcomes in which one of his opponents wins or in which nobody finds the object.

The main result of this paper is that the set of subgame perfect equilibria is exactly the set of greedy strategy profiles, also known as myopic strategy profiles, i.e. those strategy profiles in which the active player always selects one of the most likely positions containing the object. 
The key to this result is to show that by playing a greedy strategy, each player can guarantee that he wins with a probability at least as much as the probability of being active at each period.

The rest of this paper is divided as follows. In Section \ref{literature} we refer to related literature. In Section \ref{model} we introduce the model and in Section \ref{example} we present an illustrative example 
In Section \ref{greedy and SPE} we study the winning probabilities of players playing a greedy strategy, and we present a characterization of the subgame perfect equilibria. We show that the set of subgame perfect equilibria is exactly the set of greedy strategy profiles.
In Section \ref{Extensions} we discuss some extensions of the model and see to which extend the main result still holds or not.
The conclusion is in Section \ref{Conclusion}.

\subsection{Related literature} \label{literature}

The field of search problems is one of the original disciplines of Operations Research, with various applications such as military problems, R\&D problems or patent races, and many of these models involve multiple decision makers.
In the basic settings, the searcher's goal is to find a hidden object, also called the target, with maximal probability or as soon as possible. 
By now, the field of search problems has evolved into a wide range of models. 
The models in the literature differ from each other by the characteristics of the searchers and of the objects. 
Concerning objects, there might be one or several objects, mobile or not, and they might have no aim or their aim is to not be found. 
Concerning the searchers, there might be one or more. 
When there is only one searcher, the searcher faces an optimization problem. 
When there are more than one searchers, they might be cooperative or not.
If the searchers cooperate, their aim is similar to the settings with one player: they might want to minimize the expected time of search, the worst time, or some search cost function. 
If the searchers do not cooperate, the problem becomes a search game with at least two strategic non-cooperative players, and hence game theoretic solution concepts and arguments will play an important role. For an introduction to search games, we refer to  \cite{alpern2006theory}, \cite{gal1979search}, \cite{gal2010search}, \cite{gal2013search}, \cite{garnaev2012search},  and for surveys see \cite{benkoski1991survey} and  \cite{hohzaki2016search}.

The interface between Operations Research and Game Theory is attracting much attention nowadays. Examples are (PPAD) computational complexity of Nash equilibria (see \cite{papadimitriou1990graph}, \cite{papadimitriou1994complexity},
\cite{koller1996efficient}, \cite{chen20053}, \cite{savani2006hard} and \cite{chen2006settling}), price of anarchy (see \cite{roughgarden2005selfish}, \cite{bichescu2009vendor}), congestion games (see \cite{harks2012existence}, \cite{harks2010computing}, \cite{harks2016logarithmic} and \cite{harks2018uniqueness}) and manufacturer-retailer problems (see \cite{el2013dynamic} and \cite{amrouche2009shelf}).

There are many different types of search games in the literature due to variations in characteristics of the searchers and of the object. \cite{pollock1970simple} introduced a simple search game with one searcher, in which an object is moving across two locations referred to as ``state 1'' and ``state 2'' respectively, according to a discrete-time Markov chain. 
The objective of the searcher can be either to minimize the expected number of looks to find the object, or to maximize the probability of finding the object within a given horizon. The model allows for overlooking probabilities, which means that even if the searcher chooses the correct location, he may fail to find the object there.
Instead, \cite{nakai1973model} investigates the search problem with three states, while assuming perfect detection, so not accounting for overlooking probabilities.
\cite{flesch2009optimal} investigates a search game similar to \cite{pollock1970simple}, but in addition to searching for the object in state 1 or state 2, the player has a third option which is to wait. This option is costless whereas searching is costly. Waiting could induce a favourable probability distribution over the two states next period. 
They find a unique optimal strategy characterized by two thresholds.
\cite{jordan1997optimal} in his PhD thesis studies the structural properties of the optimal strategy, where the goal is to find the object while minimizing the search costs. Thereby he derived some properties of the optimal strategy for the search problem with a finite set of states in the no-overlook case and for the case where each state has the same overlooking probability and cost.
\cite{assaf1994dynamic} investigates a search problem with two states in which the object moves according to a continuous-time Markov process. 
His objective is to find the object with a minimal expected cost, where the ``real time" until the object is found is also taken into account in the cost structure.
A competitive environment with more searchers and a static object is considered in \cite{nakai1986search}.
A classical reference for an overview of search games is the book of \cite{alpern2013search}.
For a recent paper on search games, we refer to \cite{garrec2020search}. They model the search game as a zero-sum two-person stochastic game where the first player is looking for the other one. They provide upper and lower bounds on the value of the game.

\subsection{The Model} \label{model}

\noindent
\textbf{The Game.} An object is moving over a finite set $S=\{1,...,n\}$ according to a time-varying Markov chain. The initial distribution of the object is given by $\pi=(\pi_s)_{s\in S}$, and the transition probabilities at each period $t\in\N$ are given by the $S\times S$ transition matrix $P_t$, where entry $P_t(s,s')$ is the probability that the object moves to state $s'$, given it is in state $s$ at period $t$.

Let $I$ denote a set of players, who compete to find the object. We assume that $I\subseteq \N$; so the set $I$ can be either finite or countably infinite. The players do not observe the current state of the object, but they know the initial distribution $\pi$ and the transition matrices $P_t$, for each $t\in\N$. At each period $t\in \N$ one of the players is active: player $i$ is active with probability $q_i > 0$, where $\sum_{i \in I} q_i = 1$. The active player chooses a state $s_t\in S$. If the object is in state $s_t$, then the active player finds the object and wins the game. Otherwise the game enters period $t+1$. We assume that each player observes the actions chosen by his opponents, and knows the probabilities $q_i$, for each $i\in I$.

The goal of each player is to find the object and win the game.
For each player the game has three possible outcomes. 
The first possible outcome is that the player himself finds the object and wins the game. 
The second outcome is that one of his opponents finds the object and wins the game. 
The third outcome is that no one finds the object. 
In this game each player prefers the first outcome, but is indifferent between the second and third outcomes. This means that players do not have opposite interests, which makes it a non-zero-sum game.

\textbf{Actions \& Histories.} The action set for each player $i$ is $A_i=S$. Thus, a history at period $t\in \N$ is a sequence $h_t=(s_1,\ldots,s_{t-1})\in S^{t-1}$ of past actions. By $H_t=S^{t-1}$ we denote the set of all histories at period $t$. Note that $H_1$ consists of the empty sequence. Given a history $h_t$, with the knowledge of the initial probability distribution $\pi$ and the transition matrices $P_1,\ldots,P_{t-1}$, the players can calculate the probability distribution of the location of the object at period $t$. 

\textbf{Strategies.} A strategy for player $i$ is a sequence of functions $\sigma_i = (\sigma_{i,t})_{t\in \N}$ where $\sigma_{i,t} \colon H_t\rightarrow \Delta(S)$ for each period $t\in\N$. The interpretation is that, at each period $t\in \N$, if player $i$ becomes the active player, then given the history $h_t$, the strategy $\sigma_{i,t}$ recommends to search state $s\in S$ with probability $\sigma_{i,t}(h_t)(s)$. We denote by $\Sigma_i$ the set of strategies of player $i$. We say that a strategy is pure if, for any history, it places probability 1 on one action.\\
A strategy is called greedy if, for any history, it places probability 1 on the most likely states.
In the literature the greedy strategy is sometimes also called the myopic strategy.

\textbf{Winning probabilities.}
Consider a strategy profile $\sigma=(\sigma_i)_{i\in I}$. The probability under $\sigma$ that player $i$ wins is denoted by $u_i(\sigma)$. Note that \[\sum_{i \in I} u_i(\sigma)=1-\text{Prob}_{\sigma}[\text{object never found}].\]
The aim of player $i$ is to maximize $u_i(\sigma)$. If the object has not been found before period $t$, and the history is $h_t$, the continuation winning probability from period $t$ onward is denoted by $u_i(\sigma)(h_t)$ for player $i$.

\textbf{Subgame perfect equilibrium.} A strategy $\sigma_i$ for player $i$ is a best response to a profile of strategies $\sigma_{-i}$ for all other players if $u_i(\sigma) \geq u_i(\sigma'_i,\sigma_{-i})$ for every strategy $\sigma'_i \in \Sigma_i$. A strategy profile $\sigma=(\sigma_i)_{i\in I}$ is called an equilibrium of the game if $\sigma_i$ is a best response to $\sigma_{-i}$ for each player $i \in I$. 

A strategy $\sigma_i$ for player $i$ is a best response in the subgame at history $h$ to a profile of strategies $\sigma_{-i}$ for all other players, if $u_i(\sigma)(h) \geq u_i(\sigma'_i,\sigma_{-i})(h)$ for every strategy $\sigma'_i \in \Sigma_i$. A strategy profile $\sigma=(\sigma_i)_{i\in I}$ is called a subgame perfect equilibrium if, at each history $h$, in the subgame at history $h$ the strategy $\sigma_i$ is a best response to $\sigma_{-i}$  for each player $i \in I$. In other words, $\sigma$ is a subgame perfect equilibrium if it induces an equilibrium in each subgame.

\subsection{An illustrative example} \label{example}
\begin{figure}[h]
\begin{center}
\begin{tikzpicture}[->, >=stealth', auto, semithick, node distance=3cm]
\tikzstyle{every state}=[fill=white,draw=black,thick,text=black,scale=1]
\node[state]    (1)                     {$1$};
\node[state]    (2)[right of=1]   {$2$};
\path
(1) edge[bend left]          node{$1-c$}     (2)
    edge[loop left]        node{$c$}    (1)
(2) edge[loop right]           node{$0$}     (2)
    edge[bend left]      node{$1$}     (1);
\end{tikzpicture}
\end{center}
\caption{An illustrative example} 
\label{fig:SP1} 
\end{figure}
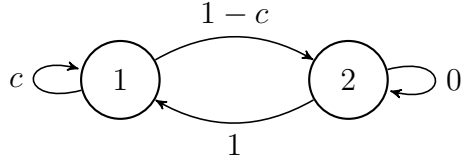

In the model we introduced the greedy strategies. However, it is not clear a priori if those strategies are relevant, nor if one greedy strategy can be better or worst than another one. 
Now we examine an example in order to try to answer those questions. 
Consider the following game with a parameter $c\in (0,1)$, with two states and two players. At each period player 1 is active with probability $q_1=0.99=1-q_2$. 
The initial probability of the location of the object is $\pi=(c,1-c)$ and the transition matrices are defined as follows: for all $t \in \N$,
\begin{equation*}
P_t=P=
\begin{pmatrix}
c & 1-c \\
1 & 0 
\end{pmatrix}.
\end{equation*}
The induced Markov chain of this game is depicted in Figure \ref{fig:SP1}.

We discuss two cases.

\textbf{Case 1.}
Consider the case in which $0<c< \tfrac{1}{3}$. At period $t=1$ it is more likely for the object to be in state 2 as $c<1-c$.

Intuitively, if player 2 gets the chance to be active at period 1, he should look in state 2 as his chance to play later is quite low. 
However, it is not immediately clear what player 1 should do at period 1.
On the one hand, if he looks in state 1 and he does not find the object, then he finds it in period 2 with probability 1 if he can play, which is likely to happen as $q_1=0.99$. On the other hand he might also simply want to maximise his chance to find the object at period 1 by looking at box 2. 

More precisely, assume first that player 2 is active at period 1. If he looks at state 1, with probability $c$ he finds the object and with probability $1-c$ the active player at period 2 can find the object by looking at state 1.
In this case, player 2 finds the object with probability $c+q_2\cdot (1-c)$.
However, if he looks at state 2, with probability $1-c$ he finds the object and with probability $c$ the active player at period 2 can find the object with probability $1-c$ by looking at state 2.
In this case, player 2 finds the object with probability at least $1-c+q_2\cdot c \cdot (1-c)> c+q_2\cdot (1-c)$ as $q_2=0.01$ and $c<\tfrac{1}{3}$. 
Thus, it is strictly better for player 2 to look at state 2 at period 1. The same holds for each subgame in which the probability distribution of the object is $(c,1-c)$. 

Now assume that player 1 is active at period 1. If he looks at state 1, then similarly to the analysis of player 2, he finds the object with probability $c+q_1\cdot (1-c)$. 
However, if player 1 decides to choose state 2 at each period he is active, by our analysis of player 2, it follows that player 2 will also do the same. Hence player 1 finds the object with probability $(1-c)+q_1 \cdot c \cdot (1-c)+q_1 \cdot c^2\cdot (1-c)+\ldots=1-c+q_1\cdot c>c+q_1\cdot (1-c)$ as $q_1<1$. So it is better for player 1 to choose state 2 at period 1. 

From our discussion, it follows that every Nash equilibrium induces the play in which at each period the active player chooses state 2. 
If a deviation occurs to state 1 and the object is not found, the next active player chooses state 1 and wins the game. This is a greedy strategy profile as at each period the most likely state is chosen.

\textbf{Case 2.} Consider the case in which $c=1/2$.

In this case, both states are equally likely, so a greedy strategy can choose any state at period 1. 
Note that at any period, if state 1 is chosen and the object is not found, a greedy strategy finds the object at the next period in state 1. 
However, if state 2 is chosen and the object is not found, a greedy strategy can choose any of the two states at the next period. 
Therefore, there are many greedy strategy profiles, and all can induce different plays. 
It is then natural to ask whether every greedy strategy profiles induce the same winning probabilities. 
As we will show later, every greedy strategy profile leads to the winning probability $q_1$ for player 1 and $q_2$ for player 2.

\section{Greedy strategies and Subgame perfect equilibria} \label{greedy and SPE}

In this section we prove our main result, which is that the set of subgame perfect equilibria is exactly the set of greedy strategy profiles. We start by introducing an intermediate result related to the winning probability guarantees of the players who play a greedy strategy.

\begin{prop} \label{winningproba}
Consider a strategy profile $\sigma=(\sigma_j)_{j\in I}$ and a player $i\in I$. If $\sigma_i$ is a greedy strategy, then under $\sigma$, the object is found with probability 1 and player $i$ wins with probability at least $q_i$. Consequently, if $\sigma$ is a greedy strategy profile, then under $\sigma$, each player $i$ wins with probability $q_i$.
\end{prop}

\textsc{Proof.} Let $\sigma=(\sigma_j)_{j\in I}$ be a strategy profile in which player $i$ plays a greedy strategy.

\textsc{Step 1.} Under $\sigma$, the object is found with probability 1.

\textsc{Proof of Step 1.} Consider an arbitrary period $t\in\N$ and suppose that the object has not been found yet. Then, player $i$ plays with probability $q_i$, and when he plays he wins with probability at least $1/n$. This implies that the object is found at period $t$ with probability at least $q_i/n$. Since this holds for each period $t$, under $\sigma$ the object is found with probability 1. 

\textsc{Step 2.} Under $\sigma$, player $i$ wins with probability at least $q_i$.

\textsc{Proof of Step 2.} Consider a period $t\in\N$ and a history $h\in H_t$. For each state $s\in S$, let $z_h(s)$ denote the probability that the object is in state $s$ at period $t$ given the history $h$. Let $z^*_h=\max_{s\in S}z_h(s)$. 

For each player $j\in I$, let $w_j(h)$ denote the probability that after history $h$ player $j$ becomes the active player and by using the mixed action $\sigma_j(h)$ he wins immediately. Since $\sigma_i$ is a greedy strategy, we have $w_i(h)=q_i\cdot z^*_h$. For each other player $j\neq i$ we have $w_j(h)\leq q_j\cdot z^*_h$. Hence, given that the object is found at period $t$ after the history $h$, the conditional probability that player $i$ finds it is 
\[\frac{w_i(h)}{\sum_{j\in I}w_j(h)}\,\geq\,\frac{q_i\cdot z^*_h}{\sum_{j\in I}q_j\cdot z^*_h}\,=\,q_i.\]

Since this holds for every period $t$ and every history $h\in H_t$, and since by Step 1 the object is found under $\sigma$ with probability 1, player $i$ wins with probability at least $q_i$. \qed

\begin{theorem} \label{thm:characterization}
The set of subgame perfect equilibria is exactly the set of greedy strategy profiles.
\end{theorem}

\textsc{Proof.} The proof is divided in two steps.

\textsc{Step 1.} Every greedy strategy profile is a subgame perfect equilibrium. 

\textsc{Proof of Step 1.} Let $\sigma$ be a greedy strategy profile and consider a subgame at a history $h\in H$. We show that $\sigma$ is an equilibrium in this subgame. Consider a player $i$ and a deviation $\sigma'_i$. By Proposition \ref{winningproba}, the strategy $\sigma_j$ guarantees to each player $j\in I$ that he wins with probability at least $q_j$, even if another player deviates. Since $\sum_{j\in I} q_j =1$, we find $u_i(\sigma)(h)=q_i$ and $u_i(\sigma'_i,\sigma_{-i})\leq q_i$. Thus, the deviation $\sigma'_i$ is not profitable.

\textsc{Step 2.} Every subgame perfect equilibrium is a greedy strategy profile. 

\textsc{Proof of Step 2.} Assume by way of contradiction that there exists a subgame perfect equilibrium $\sigma=(\sigma_i)_{i\in I}$, which is not a greedy strategy profile. 
Suppose that $\sigma_i$ is not greedy at history $h$.
Let strategy $\sigma'_i$ be a one-deviation from $\sigma_i$ at history $h$, under which player $i$ plays greedy at history $h$. So $\sigma'_{i}(h) \neq \sigma_{i}(h)$ and $\sigma'_{i}(h') = \sigma_{i}(h')$ for every $h' \in H \setminus \{h\}$.\\
Let $z_h(\sigma_i(h))$ denote the probability that player $i$ finds the object immediately with $\sigma_i$ at history $h$ given that he is the active player. 
Similarly, let $z_h(\sigma_{-i}(h))$ denote the probability that one of the opponents of player $i$ finds the object immediately with $\sigma_{-i}$ at history $h$ given that one of the opponents is active.
\\
Then we have
\begin{align*}
u_i(\sigma_i, \sigma_{-i})(h) 
&= q_i \cdot [ z_h(\sigma_i(h)) + (1-z_h(\sigma_i(h))) \cdot q_i)] + (1-q_i) \cdot (1-z_h(\sigma_{-i}(h))) \cdot q_i\\
&= q_i \cdot  z_h({\sigma}_i(h)) \cdot (1-q_i) + q_i \cdot [q_i + (1-q_i) \cdot (1-z_h(\sigma_{-i}(h)))].
\end{align*}

The first equality follows from the following argument. Player $i$ has a probability of $q_i$ of becoming the active player.
Then, he wins immediately with probability $z_h(\sigma_i (h))$ or he does not find the object immediately with probability $1-z_h(\sigma_i (h))$ and still wins in the future with probability $q_i$ by Proposition \ref{winningproba}. 
However, with probability $1-q_i$ one of his opponents is active. 
Then, player $i$ still has a chance to win.
His opponents fail to find the object immediately with probability $1-z_h(\sigma_{-i}(h))$ and then player $i$ can again win in the future with probability $q_i$. 

Similarly, for the deviation $\sigma'_i$ we have 
\begin{align*}
u_i(\sigma'_i,{\sigma}_{-i})(h) &= q_i \cdot [ z_h(\sigma'_i(h)) + (1-z_h(\sigma'_i(h))) \cdot q_i] + (1-q_i) \cdot (1-z_h(\sigma_{-i}(h)) \cdot q_i\\
&= q_i \cdot  z_h({\sigma'}_i(h)) \cdot (1-q_i) + q_i \cdot [q_i + (1-q_i) \cdot (1-z_h(\sigma_{-i}(h))].
\end{align*}

Since $z_h(\sigma'_i(h)) > z_h(\sigma_i(h))$ holds due to the fact that $\sigma'_i$ is greedy at $h$ but $\sigma_i$ is not, we find $u_i(\sigma'_i,{\sigma}_{-i})(h)>u_i(\sigma)(h)$. This however contradicts the assumption that $\sigma$ is a subgame perfect equilibrium. \qed

\section{Extensions \& Variations} \label{Extensions} 
In this section we consider several extensions of the model and discuss if the main theorem still holds or not.

\subsection{Extensions where our results still hold}

$\bullet$ \textbf{History dependent transitions.} In the model description we assumed that the object moves according to a time-varying Markov chain. A more general situation is when the transition probabilities at each period $t$ can also depend on (i) the sequence of past choices of the players, (ii) the sequence of past active players, and (iii) the sequence of states visited by the object in the past. Without any modification in the proofs, our main result, Theorem \ref{greedy and SPE}, remains valid.

$\bullet$ \textbf{Overlooking.} Our main result, Theorem \ref{greedy and SPE}, can be extended to a model with overlooking probabilities. Overlooking means that the object is in the chosen state, but the active player ``overlooks" it and therefore does not find it. Note that players cannot distinguish between overlooking and searching in a vacant state. For this extension, let $\delta_s<1$ denote the probability of overlooking the object in state $s$, for each $s\in S$. 

The overlooking probabilities need to be taken into account when defining greedy strategies. A greedy strategy should choose a state $s$ that maximizes the probability of containing the object times the probability of not overlooking the object, i.e. it should maximize the immediate probability of finding the object. 

$\bullet$ \textbf{No active player.}
The model could also be adjusted for the probability that no player is active, i.e. $\sum_{i \in I} q_i < 1$. Then, $r:= 1-\sum_{i \in I} q_i$ is the probability that no player is active at a certain period. We still assume that $q_i>0$ for each player $i$, so that $r<1$. Our main result, Theorem \ref{greedy and SPE}, would still hold, as the key properties of Proposition \ref{winningproba} remain valid.

\subsection{Variations where results break down or need adjustment}

$\bullet$ \textbf{Robustness to finite horizon and discounting.} Let $\sigma$ be a greedy strategy profile. By Theorem \ref{greedy and SPE}, $\sigma$ is a subgame perfect equilibrium, and in particular, an equilibrium. By standard arguments, for every error-term $\veps>0$, the strategy profile $\sigma$ is an $\veps$-equilibrium if the horizon of the game is finite but sufficiently long. Also, $\sigma$ is a subgame perfect $\veps$-equilibrium if, instead of maximizing the probability to find the object, now each player $i$ maximizes $\sum_{t=1}^\infty \delta^{t-1}\cdot z_{i,t}$, where $\delta\in(0,1)$ is a discount factor and $z_{i,t}$ is the probability that player $i$ finds the object at period $t$. The main reason why this robustness property holds is that, as long as at least one player plays a greedy strategy, the object will be found at an exponential rate, so essentially in finite time. 

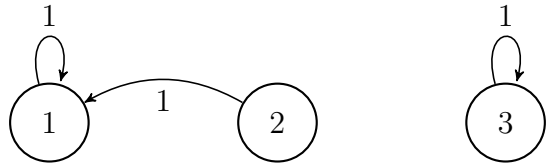
\begin{figure}[h]\begin{center}
\begin{tikzpicture}[->, >=stealth', auto, semithick, node distance=3cm]
\tikzstyle{every state}=[fill=white,draw=black,thick,text=black,scale=1]
\node[state]    (1)                     {$1$};
\node[state]    (2)[ right of=1]   {$2$};
\node[state]    (3)[ right of=2]   {$3$};
\path
(1) edge[loop above]     node{$1$}         (1)
(3) edge[loop above]                node{$1$}         (1)
(2) edge[bend right]    node{$1$}         (1);
\end{tikzpicture}
\end{center}
\caption{A 2-period game}
\label{fig:SP2} 
\end{figure}

Note however that a greedy strategy is not necessarily a 0-equilibrium on finite horizon. Indeed, consider the game represented in Figure \ref{fig:SP2} over $T=2$ periods. There are three states. The initial probability distribution of the object is given by $\pi=(1/3+\varepsilon,1/3+2\varepsilon,1/3-3\varepsilon)$, where $\varepsilon>0$ is small enough so that $0.99\cdot (2/3+3\veps)\leq 2/3$, and the transitions at period 1 are given by the arrows in Figure \ref{fig:SP2}. There are two players, and player 1 plays at each period with probability $q_1=0.99$.   

The greedy strategy profile is to choose state 2 at period 1 and state 1 at period 2. Under this strategy profile, player 1 finds the object with probability $q_1\cdot (2/3+3\veps)\leq 2/3$. However, it would be a profitable deviation for player 1 to choose state 3 at period 1 and state 1 at period 2. Indeed, player 1 is the active player at both periods with probability $(q_1)^2$, and in that case this strategy finds the object with probability 1.

$\bullet$ \textbf{Period-dependent probabilities $q_i$.} It is crucial for our proofs that the probabilities with which the players become active do not depend on the period. The reason is that in this setting, players have to balance between finding the object at the present period and putting the opponent into a difficult position in the next period.

$\bullet$ \textbf{Infinitely many states.} Assume $S=\N$. It is easy to see that it can be impossible to find the object with probability 1 if the transition law of the object is diffuse. Suppose that the object starts in states 1, 2 and 3 with equal probability, and then moves to $3^{t+1}$ new states with equal probability at each period $t$. In that case, the highest probability to find the object is at most $1/3+1/9+1/27+\cdots=1/2$. 

As a consequence, a greedy strategy profile might not be a Nash equilibrium, and vice-versa. Indeed, it might be very important to choose a certain state at period 1, even if it has a low probability of containing the object, so that such a diffusion of the object is prevented.

\section{Conclusion} \label{Conclusion}
In this paper we examined a discrete time and space search game with multiple competitive searchers who look for one object moving over finitely many locations. We showed that if each player has the same probability to play at each period, the set of subgame perfect equilibria is exactly the set of greedy strategy profiles (cf. Theorem \ref{thm:characterization}).
We discussed several variations, such as the finite truncation of the game, the discounted version of the game, cases with infinitely many states, overlooking probabilities and we examined the possibility that no player is active at a period. A challenging task would be to investigate stochastic search games when the probability of a player to be active depends on the history.


\bibliography{Searchforamovingtargetinacompetitiveenvironment.bib}{}
\bibliographystyle{plain}

\end{document}